\documentclass[a4paper]{article}
\usepackage{amsmath, amssymb, eucal}

 \newtheorem{thm}{Theorem}
 
 \newtheorem{lem}[thm]{Lemma}
 \newtheorem{prop}[thm]{Proposition}
 \newtheorem{defn}[thm]{Definition}

 \newcommand{\Complex}{\mathbb{C}}

 \newcommand{\smnoind}{{\smallskip\noindent}}
 
 \newcommand{\id}{{\rm id}}

 \newcommand{\ti}{\tilde}

\begin{document}
\title{Examples of amenable Kac system}
\author{Chi-Keung Ng}
\date{\today}
\maketitle

\begin{abstract}
By giving an interesting characterisation of amenable multiplicative 
unitaries, we show in a simple way that bicrossproducts 
of amenable locally compact groups is both amenable 
and coamenable. 
\end{abstract}

\noindent {\small 1991 AMS Mathematics Classification number:
Primary: 46L05, 46L55}

\bigskip

\bigskip

Amenability of Kac algebras has been studied in \cite{ES}.
In \cite{BS}, amenability of regular multiplicative unitaries was also 
defined. 
Recently, we studied in \cite{Ng-amen} and \cite{Ng-amen-kac} amenable 
Hopf $C^*$-algebras and amenable Kac algebras. 
However, there wasn't any non-trivial examples 
of amenable Kac algebras explicitly stated in any literature
so far.

\medskip

In this note, we will give a simple proof for the following 
expected statement (we do not know if another proof for this statement 
is already known -- note that the ``reduced algebra'' of the bicrossproduct is 
the reduced crossed product, one would imagine the ``full algebra'' to be 
the full crossed product but it seems not clear to us whether this is true): 
in the construction of the ``bicrossproduct'' (in \cite[\S 1]{BS2}) 
of two locally compact groups, 
if one of the groups is amenable, then the ``bicrossproduct'' is 
amenable (or coamenable -- depended on the convention). 
In fact, we prove this by using a simple but seeming powerful characterisation 
for amenable multiplicative unitaries (Proposition \ref{cha>amen}).

\medskip

Let us first recall from \cite[2.1 \& 2.2(b)]{Ng-morph}
as well as \cite[A.13(c)]{BS}
the following definitions. 

\medskip

\begin{defn}
\label{amen-hopf}
(a) Let $V$ be a multiplicative unitary on $H$ (in the sense of 
\cite[1.1]{BS}). 
Then $V$ is called a \mbox{\em $C^*$-multiplicative unitary} if for any 
representation $X$ and co-representation $Y$ of $V$ (see \cite[A.1]{BS}),
the sets $\hat S_X = \overline{\{(\id\otimes \omega)(X): \omega\in 
{\cal L}(H)_*\}}$ and $S_Y = \overline{\{(\omega\otimes \id)(Y): 
\omega\in {\cal L}(H)_*\}}$ are $C^*$-algebras such that 
$X\in M(\hat S_X\otimes S_V)$ and $Y\in M(\hat S_V\otimes S_Y)$. 
(Recall that in this case, $S_V$ and $\hat S_V$ are Hopf $C^*$-algebras
with coproducts $\delta$ and $\hat \delta$ defined by the formulas 
in \cite[3.8]{BS}.)

\smnoind
(b) Let $V$ be a $C^*$-multiplicative unitary. 
For any $C^*$-algebra $A$, a unitary 
$U\in M(A\otimes S_V)$ is called a {\em unitary corepresentation} if 
$(\id\otimes \delta)(U) = U_{12}U_{13}$ 
(recall that unitary corepresentations 
of $S_V$ corresponds directly to representations of $V$). 
Moreover, there exists an 
{\em universal object $(\hat S_p, V')$ for unitary corepresentations
of $S$} in the sense 
that $\hat S_p$ 
is a Hopf $C^*$-algebra and $V'\in M(\hat S_p\otimes S_V)$  is a unitary 
corepresentation of $S_V$ such that 
for any unitary corepresentation $X$ of $S_V$, there exists a canonical 
surjective $*$-homomorphism $\pi_X$ from $\hat S_p$ to 
$\hat S_X$ such that $X = (\pi_X\otimes \id)(V')$. 

\smnoind
(c) A $C^*$-multiplicative unitary $V$ is said to be {\em amenable} if 
the trivial representation $\pi_1$ is weakly contained in $\pi_V$.
On the other hand, $V$ is called {\em coamenable} if $\Sigma V^* \Sigma$ is 
amenable. 
\end{defn}

\medskip

Any regular or semi-regular and balanced or manageable multiplicative 
unitary is a $C^*$-multiplicative unitary. 
Note that part (c) is only one of the many equivalent formulations for 
amenability (see e.g. \cite{Ng-amen} for some of them or \cite{Ng-amen-kac} 
for more of them in the case of Kac algebras). 

\medskip

To prove the main theorem, we need Proposition \ref{cha>amen} 
which is a partial generalisation of 
\cite[4.8(XII)]{Ng-amen-kac} (in the sense that we only concern with 
one dimensional representations instead of finite dimensional ones). 
Note that this proposition implies that the amenability 
of a $C^*$-multiplicative unitary $V$ depends solely on the $C^*$-algebraic 
structure of $\hat S_V$. 
The proof of the proposition require the following lemma 
which itself follows from \cite[5.5]{Ng-dual} (note that 
if $L$ is the canonical representation of $S$ on ${\cal L}(H)$, 
then $(L,\pi_V)$ is a $V'$-covariant representation 
in the sense of \cite[4.1]{Ng-dual}). 

\medskip

\begin{lem}
\label{support}
Let $V$ be a $C^*$-multiplicative unitary on a Hilbert space $H$. 
Suppose that $\ti\pi_{V}$ is the extension of $\pi_{V}$ to 
$\hat S_p^{**}$ (considered as a representation 
on ${\cal L}(H)$) and $e$ is the support projection of $\ti\pi_{V}$. 
Then $\hat \delta_p^{**}(e) \geq e\otimes 1$
(where $\hat \delta_p$ is the coproduct on $\hat S_p$).
\end{lem}

\begin{prop}
\label{cha>amen}
Suppose that $V$ is a $C^*$-multiplicative unitary on $H$.  
Then $V$ is amenable if and only if $\hat S_V$ has 
an one dimensional representation. 
\end{prop}
\noindent
{\bf Proof:}
Since the amenability of $V$ will imply the existence of a counit for 
$\hat S_V$ (see e.g. \cite[3.5]{Ng-amen}), we need only to show 
the converse. 
Suppose that $\psi$ is an one dimensional representation of $\hat S_V$. 
Then $U=(\psi\otimes \id)(V)\in M(S_V)$ is a 
unitary corepresentation of $S_V$. 
This means that $U$ is a group-like element in $M(S_V)$ and 
$\pi_U = \psi \circ \pi_V$. 
Moreover, $U^*$ is again a group-like unitary and hence a unitary 
corepresentation of $S_V$ on $\Complex$ (or equivalently, 
a representation of $V$ on $\Complex$).  
Now $(\psi\circ \pi_V\otimes \pi_{U^*})\circ \hat \delta_p = 
(\pi_U\otimes \pi_{U^*})\circ \hat\delta_p = \pi_{UU^*} = \pi_1$. 
Using the argument of \cite[3.6]{Ng-amen-kac}, we see that 
$\ker \pi_V \subseteq \ker(\pi_V\otimes\id)\circ\hat\delta_p$. 
For the benefit of the readers, we repeat the simple argument here. 
Note that $\ker\pi_{V} = (1-e)\hat S_p^{**}\cap \hat S_p$  
and so $x=(1-e)x$ for any $x\in \ker \pi_{V}$. 
This, together with $(e\otimes 1)\hat \delta_p^{**}(1-e) = 0$ 
(given by Lemma \ref{support}), implies that 
\begin{eqnarray*}
(\pi_{V}\otimes \id)\hat \delta_p(x)
& = & (\ti\pi_{V}\otimes \id)\hat \delta_p^{**}(x)\\
& = & (\ti\pi_{V}\otimes \id)\hat \delta_p^{**}((1-e)x)
\quad = \quad (\ti\pi_{V}\otimes \id)((e\otimes 1)\hat \delta_p^{**}((1-e)x)) 
\quad = \quad 0.
\end{eqnarray*}
Hence $\ker\pi_V \subseteq \ker (\psi\otimes \pi_{U^*})\circ 
(\pi_V\otimes\id)\circ\hat\delta_p = \ker \pi_1$. 
In other words, $\pi_1$ is weakly contained in $\pi_V$. 

\bigskip

In the case of Kac algebras, we can obtain this proposition by an 
even simpler proof (without the need of Lemma \ref{support} nor 
the argument of \cite[3.6]{Ng-amen-kac})
because $\pi_{U^*} = \psi\circ\hat\kappa\circ\pi_V$. 

\medskip

From now on, $G$ is a locally compact group and $G_1$ and $G_2$ are two 
closed subgroups of $G$ such that the map $\varphi:G_1\times G_2 
\rightarrow G$ given by $\varphi(r,g)=rg$ is a homeomorphism 
from $G_1\times G_2$ onto an open dense subset $\Omega$ of $G$ 
(see \cite[\S 1]{BS2}). 
Consider the canonical actions $\alpha$ and $\beta$ of $G_1$ and 
$G_2$ on $G_1\!\!\setminus\! G$ and $G/G_2$ respectively. 
Note that both $\alpha$ and $\beta$ acts by automorphisms. 
Let $V$ be the semi-regular irreducible 
multiplicative unitary as in \cite[1.1]{BS2}. 

\medskip

\begin{thm}
\label{amen-bicross}
With the notation as in the above paragraph, if $G_1$ 
(respectively, $G_2$) is amenable, then $V$ is amenable 
(respectively, coamenable). 
\end{thm}
\noindent
{\bf Proof:}
By \cite[1.1(c)]{BS2}, we know that $\hat S_V$ is isomorphic to the reduced
crossed product $C_0(G_1\!\!\setminus\! G)\times_{\alpha,r} G_1$. 
Therefore, using Proposition \ref{cha>amen}, it suffices to show that 
$C_0(G_1\!\!\setminus\! G)\times_{\alpha,r} G_1$ has an one dimensional representation. 
Let $\phi$ be the one dimensional representation 
on $C_0(G_1\!\!\setminus\! G)$ given by the 
evaluation at the coset $[e] = G_1$. 
Let $\pi$ be trivial representation of $G_1$ (in fact, 
we can take any one dimensional representation). 
We note that $(\phi, \pi)$ is a covariant representation 
for $\alpha$. 
Indeed, for any $f\in C_0(G_1\!\!\setminus\! G)$, we have 
$\phi(\alpha_s(f)) = f(\alpha_{s^{-1}}([e])) = f([e]) = \phi(f) = 
\pi(s)\phi(f)\pi(s)^*$. 
Hence, we obtained an one dimensional representation for the full crossed product
$C_0(G_1\!\!\setminus\! G)\times_{\alpha,{\rm max}} G_1$ which is 
identical with the reduced crossed product 
$C_0(G_1\!\!\setminus\! G)\times_{\alpha,{\rm r}} G_1$
(as $G_1$ is amenable). 

\bigskip

Note that $\Omega = G$ if and only if $G_1\cong G/G_2$ and 
$G_2 \cong G_1\!\!\setminus\! G$ under the canonical maps. 
In this case, the actions $\alpha$ and $\beta$ are given by 
the commutation rule: $gr = \beta_g(r) \alpha_r(g)$ (that comes from the 
bijectivity of $\varphi$; see \cite[\S 8]{BS}). 
Moreover, $\alpha$ and $\beta$ satisfy the conditions for matching 
pairs in \cite{Maj} (except those concerning the Radon-Nikodym derivatives). 
Recall that in this case, $G$ is the bicrossproduct group as 
defined in \cite[\S 2]{Yam}. 
Note that the bicrossproduct group can only be defined for a
matching pair with $\alpha_e$ and $\beta_e$ being the identity maps 
(the formula for the inverse in \cite[pp.255]{Yam} 
requires such conditions) or equivalently,
$\alpha$ and $\beta$ act by automorphisms.

\medskip

On the other hand, in this case, 
$(L^2(G),V,U)$ (where $U$ is the unitary as in \cite[1.1(b)]{BS2})
is the Kac system bicrossproduct 
of $(L^2(G_1), X, u)$ (the canonical Kac system of $G_1$ as 
given in \cite[6.11(b)]{BS}) and 
$(L^2(G_2), Y, v)$ relative to $Z$ (see \cite[8.13, 8.15 \& 8.21]{BS}). 
Therefore, Theorem \ref{amen-bicross} applies to Kac system 
bicrossproducts of locally compact groups (as defined in \cite[\S 8]{BS}). 
If, in addition, $\alpha$ and $\beta$ in the above paragraph satisfy 
the conditions for the
Radon-Nikodym derivatives in \cite{Maj}, this Kac system actually comes 
from a Kac algebra which is also amenable and coamenable.

\noindent
Department of Pure Mathematics,
The Queen's University of Belfast,
Belfast BT7 1NN,
United Kingdom.

\smnoind
E-mail address: c.k.ng@qub.ac.uk

\end{document}